\documentclass[12pt, a4paper]{article}
\usepackage[left=1in,right=1in,top=1.0in,bottom=1in]{geometry}

\usepackage{amsmath,amssymb,amsthm,amsfonts,afterpage,hyperref}
\usepackage{amscd,subeqnarray}

\usepackage{graphicx}
\usepackage{caption}
\usepackage{subcaption}
\usepackage{authblk}
\usepackage{tikz}

\usepackage{array}
\usepackage{wrapfig}
\usepackage{color}
\usepackage{ragged2e}
\usepackage{stmaryrd}
\usepackage{siunitx}
\usepackage{multirow}
\usepackage{xcolor,cancel}
\usepackage{boldfonts}
\usepackage{symbols}
\usepackage{footmisc}
\usepackage{float}
\usepackage{setspace}
\usepackage{bm}
\usepackage{booktabs}
\usepackage{tikz, xcolor}
\usepackage{soul}
\usepackage{pdfpages} 

\usepackage[T1]{fontenc}
\usepackage[utf8]{inputenc}

\usepackage[linesnumbered,ruled,vlined]{algorithm2e}
\SetKwInput{KwInput}{Input}                
\SetKwInput{KwOutput}{Output} 

\usepackage[sort, numbers]{natbib}
\setcitestyle{square}

\usetikzlibrary{shapes,arrows, positioning,calc}	
\usetikzlibrary{arrows,snakes,backgrounds}

\newcommand{\bfa}[1]{\boldsymbol{#1}} 			
			
\newcommand{\bfeps}{\boldsymbol{\epsilon}}

\newcommand{\tr}{\text{tr}}       				%

\DeclareMathAlphabet{\mathpzc}{OT1}{pzc}{m}{it}

\newcommand{\bfn}{\boldsymbol{n}}	
	
\newcommand{\bfu}{\boldsymbol{u}}	
	
\newcommand{\bfv}{\boldsymbol{v}}

\newcommand{\bfx}{\boldsymbol{x}}

\newcommand{\bfT}{\boldsymbol{T}}		
\newcommand{\bfI}{\boldsymbol{I}}	 
\newcommand{\bfzero}{\boldsymbol{0}}

\newcommand{\bff}{\boldsymbol{f}}	
\newcommand{\bfg}{\boldsymbol{g}}	

\newtheorem{remark}{Remark}

\newtheorem*{cf}{Strong formulation}
\newtheorem*{cwf}{Continuous weak formulation}

\providecommand{\keywords}[1]
{
  \small	
  {\textit{Keywords---}} #1
}

\title{Crack-tip field characterization in nonlinearly constituted and geometrically linear elastoporous solid containing a star-shaped crack: A finite element study}

\author[1]{S. M. Mallikarjunaiah}
\author[2]{Kun Gou}

\affil[1]{Department of Mathematics \& Statistics\\
Texas A\&M University-Corpus Christi\\
Corpus Christi, TX 78412, USA\\
Email: m.muddamallappa@tamucc.edu}
\affil[2]{Department of Computational, Engineering, and Mathematical Sciences\\
Texas A\&M University-San Antonio\\
San Antonio, TX 78224, USA\\
Email: kun.gou@tamusa.edu}
\date{}

\begin{document}

\maketitle

\begin{abstract}
This paper introduces a three-dimensional (3-D) mathematical and computational framework for the characterization of crack-tip fields in star-shaped cracks within porous elastic solids. A core emphasis of this model is its direct integration of density-dependent elastic moduli, offering a more physically realistic representation of engineering materials where intrinsic porosity and density profoundly influence mechanical behavior. The governing boundary value problem, formulated for the static equilibrium of a 3-D, homogeneous, and isotropic material, manifests as a system of second-order, quasilinear partial differential equations. This system is meticulously coupled with classical traction-free boundary conditions imposed at the complex crack surface. For the robust numerical solution of this intricate nonlinear problem, we employ a continuous trilinear Galerkin-type finite element discretization. The inherent strong nonlinearities arising within the discrete system are effectively managed through a powerful and stable {Picard-type linearization scheme}. The proposed model demonstrates a remarkable ability to accurately describe the full stress and strain states in a diverse range of materials, crucially recovering the well-established classical singularities observed in linearized elastic fracture mechanics. A comprehensive numerical examination of tensile stress, strain, and strain energy density fields consistently reveals that these quantities attain their peak values in the immediate vicinity of the crack tip, an observation that remarkably aligns with established findings in standard linearized elastic fracture mechanics. Furthermore, the robust framework presented in this article provides a powerful tool for predicting material failure, allowing for the application of conventional local fracture criteria to rigorously investigate both the quasi-static and dynamic evolution of crack tips, thereby advancing our understanding of fracture initiation and propagation in complex porous media.
\end{abstract}

\noindent \keywords{Density-dependent moduli, Preferential stiffness, Implicit constitutive relation, Star-shaped crack, Porous elastic solid, Finite element method}

\section{Introduction}\label{sec:Intro}
Understanding the behavior of star-shaped cracks is of paramount importance in the fields of fracture mechanics and material science. These complex crack geometries frequently emerge from challenging loading scenarios or inherent material imperfections, as observed in impact events on brittle substances like glass windshields or during the manufacturing of materials such as steel slabs \cite{theocaris1977new}. Accurately predicting material failure hinges on comprehending the unique and intricate stress concentrations that develop at the multiple tips of such cracks. Unlike simpler linear cracks, the interaction among the various crack branches in a star configuration generates a far more complex stress field, offering a more realistic representation of certain failure mechanisms. Furthermore, the analysis of star-shaped cracks provides critical insights into crack propagation and branching behavior. Investigating how these cracks initiate, advance, and potentially merge or interact is fundamental for evaluating the remaining structural integrity of components \cite{vandenberghe2013star,andersson1969stress}. The knowledge gained from studying star-shaped cracks can directly inform the development of more accurate predictive models for crack growth, which are essential for designing durable and safe materials and structures. Ultimately, the capability to precisely model and forecast the response of star-shaped cracks significantly advances both fracture mechanics theories and computational methodologies. The inherent complexity of these geometries often necessitates the application of sophisticated numerical techniques, such as finite element analysis. Research in this area continually expands our understanding of material behavior under extreme conditions, leading to the creation of more robust design principles and enhanced failure analysis tools applicable across a broad spectrum of engineering disciplines.

In the traditional framework of elasticity theory, the majority of engineering materials are conventionally modeled as compressible. However, a specific class of materials, characterized by a {Poisson's ratio equal to $0.5$}, exhibits an intriguing property: they undergo no volume changes whatsoever under any applied deformation. Such materials, possessing an invariant density, are thus rigorously approximated as incompressible. Yet, this simplification falls short for a vast array of truly compressible materials, notably {polymers}, where density exerts a profound influence on fundamental material properties such as strength, hardness, and brittleness. Consequently, to accurately capture the elastic response of complex polymer networks and their behavior under mechanical loading, it becomes imperative to account for the {dependence of material moduli on density} \cite{rubinstein2002elasticity}.

A direct consequence of these observations is the clear necessity for a {correction to the classical elasticity model} to adequately describe the bulk behavior of numerous porous materials. While the response of many polymeric materials undergoing infinitesimal deformation has historically been modeled using the classical linearized elastic constitutive relationship (wherein Young's modulus, Poisson's ratio, and Lam\'e constants are assumed to be invariant), experimental evidence often contradicts this assumption. For instance, a study on the strength of silica gels \cite{leventis2002nanoengineering} explicitly reported a {nonlinear relationship between load and strain} during three-point bending tests on aerogels (see Figure~3 in \cite{leventis2002nanoengineering}), leading to a calculated Young's modulus that varies with conditions. Further research, as highlighted in \cite{chandrasekaran2017carbon,anjos2010tensile}, consistently demonstrates that {material moduli are indeed density-dependent}. This intricate relationship is particularly relevant in most solids where macro-level porosity is intrinsically linked to local density, with porosity and density exhibiting an inverse proportionality \cite{ALAGAPPAN2023100162}. Recognizing these limitations, Rajagopal, in a series of foundational articles including \cite{rajagopal2003,rajagopal2007elasticity,rajagopal2014nonlinear,rajagopal2007response,rajagopal2009class}, significantly {generalized the framework of elastic bodies}. Through rigorous linearization procedures, he demonstrated the derivation of {constitutive relations that are implicitly related to the linearized strain and Cauchy stress}. This implicit structure is particularly well-suited for describing the complex response of elastic bodies that are unable to convert mechanical work into thermal energy. A key advantage of these generalizations is their broad applicability; they have been shown to encompass both {Cauchy and Green elastic bodies as specific subclasses} within this broader array of elastic material models \cite{Mallikarjunaiah2015,MalliPhD2015,ortiz2014numerical,ortiz2012}.

Crucially, the aforementioned experimentally observed nonlinear material behavior, particularly the density dependence, can be effectively captured and modeled within the comprehensive framework of the {implicit theory of elasticity} \cite{rajagopal2014nonlinear,muliana2018determining,kowalczyk2019finite}. This theoretical paradigm offers the flexibility required to describe such intricate responses. For instance, a particular subclass of models derived from the general nonlinear relationships presented in \cite{rajagopal2003,rajagopal2007elasticity} has proven effective in modeling materials like {rocks} \cite{bustamante2020novel} and various {rubber-like materials} \cite{bustamante2021new}. A notable example is the constitutive model for rubber proposed in \cite{bustamante2021new}, which utilizes principal stress as the primary variable. This approach has demonstrated a superior correlation with experimental data compared to more traditional models such as Ogden's model \cite{ogden1972large}, further underscoring the power and versatility of this generalized implicit framework. A significant advancement in modeling elastic materials like rocks and concrete comes from a specific subclass of nonlinear models originally presented in \cite{rajagopal2003,rajagopal2007elasticity}. For these materials, the dependence of material moduli on density is rigorously established in more recent works such as \cite{rajagopal2021b,rajagopal2022implicit}. A key characteristic of these models is that both the linearized strain and Cauchy stress appear linearly within the constitutive relations \cite{itou2022investigation,itou2021implicit,itou2023generalization,pruuvsa2022pure,pruuvsa2023mechanical}.

Despite this linear appearance, the overall material response described by the models in \cite{rajagopal2021b,rajagopal2022implicit} is fundamentally {nonlinear}. This nonlinearity arises because the constitutive relations incorporate a bilateral product of the Cauchy stress and linearized strain. Nevertheless, the inherent structure of these models is particularly powerful: it allows for the description of a linear homogeneous material that possesses intrinsic properties, such as stiffness, which are explicitly dependent on its local density. Further exploring this domain, \cite{yoon2024finite} introduced a two-dimensional model specifically for materials whose properties are density-dependent. This study delved into the intricacies of {crack-tip fields} and {mechanical stiffness} for anisotropic materials under various loading conditions. It successfully captured a diverse range of phenomena for both tensile and in-plane shear loading. Crucially, this work reported several compelling failure scenarios that are impossible to predict using classical linear elasticity models, underscoring the necessity and advantage of these advanced nonlinear frameworks.

In this paper, we develop a comprehensive model for a three-dimensional (3-D) plate featuring an internal star-shaped crack with small opening angles. Within the framework of {nonlinear constitutive relations}, where both Cauchy stress and strain are considered, we formulate a detailed boundary value problem (BVP) tailored to this star-shaped crack geometry. To solve the resulting quasi-linear elliptic partial differential equation system, we propose a robust and convergent numerical method based on continuous Galerkin-type finite elements. This approach allows us to accurately simulate the complex mechanical behavior near the crack. Crucially, our study includes several postprocessing analyses to fully characterize the crack-tip fields. We present detailed computational results for the {crack-tip stress, strain, and strain energy density} along a line extending towards the star-crack tip. Ultimately, our work lays a vital foundation for the future development of a sophisticated fracture theory applicable to brittle materials with density-dependent moduli. These materials exhibit properties that are inherently sensitive to their local density, which is critical for accurately capturing their response to mechanical loading. 

\section{Implicit constitutive relations}\label{sec:icr}
The physical domain of the material under consideration is represented by a closed and bounded {Lipschitz domain} $\Lambda \subset \mathbb{R}^3$ in its reference configuration. The boundary of this domain, denoted as $\partial \Lambda$, is Lipschitz continuous and is rigorously partitioned into two disjoint subsets: $\partial \Lambda = \overline{\Gamma_{D}} \cup \overline{\Gamma_{N}}$. Here, $\Gamma_{D}$ signifies a non-empty Dirichlet boundary where essential boundary conditions (e.g., prescribed displacements) are applied, while $\Gamma_{N}$ represents a Neumann boundary where natural boundary conditions (e.g., prescribed tractions) are imposed. We define $\bfa{X}$ as a typical point in the reference configuration and $\bfa{x}$ as its corresponding point in the deformed configuration. The mechanical response of the material is described by the displacement field $\bfa{u} \colon \Lambda \to \mathbb{R}^3$, which relates the reference and deformed configurations through the relationship $\bfa{u} = \bfa{x} - \bfa{X}$.

For the remainder of this paper, we adopt standard notations and definitions from functional analysis, particularly those pertaining to Lebesgue and Sobolev spaces, consistent with established literature such as \cite{ciarlet2002finite,evans1998partial}. The Lebesgue space $L^{p}(\Lambda)$ for $p \in [1, \infty)$ comprises functions whose $p$-th power of the absolute value is integrable over $\Lambda$. A particularly relevant instance is $L^{2}(\Lambda)$, which denotes the space of square-integrable functions defined on $\Lambda$. This space is equipped with the inner product $\left( v, w \right) := \int_{\Lambda} v \, w \, d\bfa{x}$ and its induced norm $\| v \|_{L^2(\Lambda)} := \left(v, v \right)^{1/2}$. Furthermore, $C^{m}(\Lambda)$ represents the space of $m$-times continuously differentiable functions on $\Lambda$, where $m \in \mathbb{N}_0$. The standard Sobolev space $H^{1}(\Lambda)$ is formally defined as:
\begin{equation}
H^{1}(\Lambda) := W^{1, 2}(\Lambda) = \left\{ v \in L^{2}(\Lambda) \; : \; \partial_j v \in L^{2}(\Lambda) \text{ for } j \in \{ 1, 2, 3 \} \right\},
\end{equation}
and is endowed with its canonical inner product and induced norm. The subspace $H_0^{1}(\Lambda)$ is defined as the closure of $C_0^{\infty}(\Lambda)$ (the space of infinitely differentiable functions with compact support within $\Lambda$) under the $H^1(\Lambda)$ norm, i.e., $H_0^{1}(\Lambda) = \overline{C_0^{\infty}(\Lambda)}^{\| \cdot \|_{H^1}}$.

In this study, our primary focus is on characterizing the behavior of \textit{elastic (nondissipative) porous solids}, specifically those materials where the constitutive moduli exhibit a direct dependence on the material density. The mechanical response of such complex materials under applied loading is most accurately captured through \textit{implicit constitutive relations}, a framework pioneered by Rajagopal \cite{rajagopal2021a, rajagopal2021b,rajagopal2022implicit,rajagopal2021note} and further elaborated in his extensive body of work and the references cited therein.

A generalization of the classical elastic body, as initially defined by Rajagopal \cite{rajagopal2003, rajagopal2007elasticity} through an implicit type constitutive relation, establishes a fundamental relationship between the Cauchy stress ($\bfa{T}$) and the deformation gradient ($\bfa{F}$). This relationship is expressed in the form:
\begin{equation}
\bfa{f}(\rho, \, \bfa{T},\,\bfa{F}, \, \bfa{X}) =0,
\end{equation}
where $\bfa{f}$ denotes a tensor-valued function that implicitly defines the material's response, and $\rho$ represents the material density in the deformed configuration, while $\bfa{X}$ is the position vector in the reference configuration. For the specific case of {isotropic bodies}, this implicit response relation simplifies considerably, reducing to the form:
\begin{equation}\label{eq:implicit_relation}
\bfa{f}(\rho, \, \bfa{T},\,\bfa{B}) =0.
\end{equation}
Here, $\bfa{B}$ is the left Cauchy-Green deformation tensor, an objective measure of deformation. This implicit framework provides a powerful tool for describing materials where the stress-strain relationship cannot be explicitly inverted, offering enhanced versatility for complex material behaviors such as those observed in porous media with evolving densities.

Within the broader class of implicit constitutive relations, as generically represented by Equation \eqref{eq:implicit_relation}, a particularly relevant subclass of models describes materials where the constitutive response exhibits {linearity in both the stress tensor ($\bfT$) and the infinitesimal strain tensor ($\bfeps$)}. One such specialized constitutive model, further explored in \cite{itou2021implicit}, is given by:
\begin{equation} \label{spe_model1}
\left( 1 + \delta_3 \, \text{tr} \, \bfT \right) \; \bfeps = C_1 \, \left( 1 + \delta_1 \, \text{tr} \, \bfeps \right) \; \bfT + C_2 \, \left( 1 + \delta_2 \, \text{tr} \, \bfeps \right) \; \left( \text{tr} \bfT \right)\; \boldsymbol{I}.
\end{equation}
In this formulation \eqref{spe_model1}, the material parameters $\delta_1$, $\delta_2$, $\delta_3$, $C_1$, and $C_2$ are defined as constants. This model is notable for its inherent linearity with respect to both stress and strain, making it well-suited for describing the mechanical behavior of certain porous solids. The {classical linear elastic constitutive model} can be precisely recovered from Equation \eqref{spe_model1} by setting the parameters $\delta_1$, $\delta_2$, and $\delta_3$ to zero. In this limiting case, the constants $C_1$ and $C_2$ can be identified with the material's elastic properties as:
\begin{equation}\label{eq:material_coeff}
C_1 = \dfrac{1 + \nu}{E} = \dfrac{1}{2 \, \mu} >0, \quad C_2 = - \dfrac{ \nu}{E} <0,
\end{equation}
where $E$ represents the Young's modulus and $\nu$ is the Poisson's ratio. These elastic moduli are, in turn, intrinsically related to the Lamé constants, $\lambda$ and $\mu$, by the following well-known expressions:
\begin{equation}\label{eq:linear_lame}
\lambda = \dfrac{ E \, \nu}{(1 + \nu) (1 - 2 \nu)}, \quad \mu = \dfrac{ E }{2 (1 + \nu)}.
\end{equation}
Furthermore, by imposing specific constraints on the material parameters, such as setting $\delta_3 = 0$ and $\delta_1 = \delta_2 = \beta$ in Equation \eqref{spe_model1}, we can derive a simplified, yet significant, constitutive relation:
\begin{equation} \label{spe_model2}
\bfeps = C_1 \, \left( 1 + \delta \, \text{tr} \, \bfeps \right) \, \bfT + C_2 \, \left( 1 + \beta \, \text{tr} \, \bfeps \right) \; \left( \text{tr} \bfT \right)\; \boldsymbol{I}.
\end{equation}
Utilizing the principle of {balance of mass} in the reference configuration, the aforementioned relation can be equivalently expressed in terms of the material density $\rho$ as:
\begin{equation} \label{spe_model3}
\bfeps = \Phi_1\left( C_1, \; \rho, \; \text{tr} \, \bfeps \right) \;  \bfT + \Phi_2\left(C_2, \; \rho, \; \text{tr} \, \bfeps\right) \; \left( \text{tr} \bfT \right)\; \boldsymbol{I},
\end{equation}
with both functions $\Phi_i\left(C_i, \; \rho, \; \text{tr} \, \bfeps \right), \;\; i=1, \, 2,$ being linear in the density variable $\rho$, hence linear in $\text{tr} \, \bfeps$. The expression $\left( 1 + \delta \, \tr \, \bfeps \right)$ is a {density-dependent modulus}. This is because the trace of the strain tensor, $\tr \, \bfeps$, directly relates to the density changes between two different configurations of the material by the \textit{mass balance},
\begin{equation}\label{mb}
\rho_0 = \rho \, \left( 1 + \tr(\bfeps)  \right),
\end{equation}
where $\rho_0$ and $\rho$ are the densities of two configurations. Consequently, model \eqref{spe_model3} is particularly well-suited for describing the behavior of {porous materials} where density variations play a significant role in their mechanical response. In such materials, it's common for the {material moduli} (which describe the material's stiffness and resistance to deformation) to vary with changes in density. This stands in contrast to the {linearized elastic description of solids}, where material moduli are typically assumed to be constant and {density-independent}.

In the upcoming section, we'll delve into the development of a mathematical model specifically designed to describe the stress and strain within porous elastic materials. This model will leverage novel constitutive relations that we've established. Following this, we'll introduce a stable discretization method for the resulting system of partial differential equations. Finally, we'll highlight the key distinctions between the numerical outcomes generated by the traditional linearized elasticity model and our novel implicit model, showcasing the advantages of our new approach.

\section{Formulating the BVP for 3-D porous elastic solids}\label{sec:BVP}
This section is dedicated to establishing a comprehensive BVP that governs the mechanical behavior of a 3-D porous elastic solid. Our formulation is grounded in the unique constitutive relationship introduced previously in Section \ref{sec:icr}. A primary objective here is to meticulously characterize the intricate crack-tip stress and strain fields within the framework of these new constitutive relations. This novel relationship is distinctive in that both stress and linearized strain appear linearly, yet the material's moduli are intrinsically linked to its density. For the purposes of this analysis, we make the simplifying assumptions that the bulk material under consideration is homogeneous, isotropic, and initially in an unstrained and unstressed state. To precisely describe the stress-strain state within this new class of material models, we adopt a static problem formulation. Consequently, the fundamental principle of the balance of linear momentum, in the absence of any external body forces, simplifies to:
\begin{equation}
    - \nabla \cdot \bfT = \boldsymbol{0} \quad \mbox{with} \quad  \bfT=\frac{\mathbb{E}[\bfeps]}{1 + \beta \, \tr(\bfeps)} \text{in} \quad \Omega. \label{eq_blm}
\end{equation}
Our approach to accurately representing the material's stress response, particularly with its moduli influenced by density, involves utilizing the constitutive relation given in equation \eqref{spe_model1}. In Equation~(\ref{eq_blm}), the symbol $\mathbb{E}[\bfeps]$ is the fourth-order elasticity tensor with point-wise varying Lam\'e coefficients: 
\begin{equation}\label{eq:c1_c2_dd_model}
\lambda := \dfrac{\nu \, E}{(1 + \nu)(1-2\nu)(1+ \beta \, \tr(\bfeps))}, \quad \mu := \dfrac{E}{2(1 + \nu)(1+ \beta \, \tr(\bfeps))}.
\end{equation}
\begin{remark}
It's crucial to recognize the inherent {nonlinear dependence on strain} within the stress function defined in equation \eqref{eq_blm}. This nonlinearity stems not only from the explicit appearance of the strain function in the denominator but also from the \textit{generalized Lamé coefficients} themselves, which are implicitly strain-dependent. This contrasts sharply with traditional homogeneous elastic models. Here, the Lamé coefficients ($\lambda$ and $\mu$) are {pointwise varying}, meaning their values can differ at every location within the material. This feature is essential for accurately representing {heterogeneous materials} like composites, functionally graded materials, or geological formations, where material properties such as composition, density, or microstructure naturally vary with position. Such spatial variations fundamentally influence the stiffness and bulk modulus, leading to complex stress and strain distributions that cannot be adequately described by classical homogeneous elasticity.
\end{remark}
\begin{remark}
It's essential to recognize that while the general class of constitutive relations given by \eqref{spe_model1} isn't universally invertible, a specific selection of parameters, as used in \eqref{spe_model1}, yields an invertible material response. This invertibility is crucial as it allows us to explicitly define stress as a nonlinear function of strain, as shown in \eqref{eq_blm}. A significant area for future research involves rigorously identifying the conditions for {strong ellipticity} and {convexity} of the relations in \eqref{spe_model1} that lead to models of the type presented in \eqref{eq_blm}. Such an investigation is vital for ensuring the stability of numerical approximations when solving BVPs formulated within the framework of \eqref{spe_model1}.
\end{remark}

Consider a porous material $\Omega \subset \mathbb{R}^3$ containing the star-shaped crack under tensile loading. We seek a displacement field $\bfa{u}(\mathbf{x}) \colon \Omega \to \mathbb{R}^3$ that satisfies the following BVP:
\begin{cf}
Find $\bfa{u} \in C^2(\Omega) \cap C^0(\bar{\Omega})$ such that:
\begin{align}
-\nabla \cdot \left[ \frac{\mathbb{E}[\bfeps]}{1 + \beta \, \tr(\bfeps)}          \right] &=\bfa{f} \quad \text{in } \Omega \label{eq:governing}\\
 \bfu &= \bfu^0, \quad  \mbox{on} \quad \Gamma_D, \;\; \mbox{and}  \label{eq:dirichlet} \\
 \bfT  \bfn &= \bfg, \;\;  \mbox{on} \;\; \Gamma_N, \label{eq:neumann}\\
\bfT  \bfn &= \bfa{0}, \;\;  \mbox{on} \;\; \Gamma_C,
\end{align}
\end{cf}
Here, $\Omega$ is the open, bounded domain with boundary $\partial\Omega = \Gamma_D \cup \Gamma_N$ and $\Gamma_D \cap \Gamma_N = \emptyset$, $f(\mathbf{x})$ represents the volumetric source term, $\mathbf{n}$ is the outward unit normal vector on the boundary, $\bfa{u}^0$ specifies the prescribed dispalcement on the Dirichlet boundary $\Gamma_D$, $\bfa{g}$ specifies the prescribed flux on the Neumann boundary $\Gamma_N$, and $\Gamma_C$ represents the surfaces of the star-shaped crack. 

As \eqref{eq:governing} is inherently {nonlinear}, its direct solution through conventional analytical or numerical methods proves challenging. To overcome this, in the following section, we propose a \textit{Picard's iteration} for linearization at the differential equation level. This will be subsequently discretized using a {bilinear finite element method} to obtain a convergent numerical solution. 

\subsection{Finite element discretization of the nonlienar crack model}
\label{sec:fe_discretization} 

This section details the {finite element discretization} employed to investigate the mechanical response of a 3-D porous elastic solid, where the material properties are contingent upon its density. Our primary objective is to meticulously compare the model's predictions, as proposed herein, against those derived from the established classical linearized elasticity theory. The behavior of this porous elastic solid is governed by the BVP formulated in the preceding section, as supported by recent works such as \cite{murru2021stress,itou2021implicit,rajagopal2021b,yoon2024finite}. We proceed under the assumption that a unique strong solution exists for the governing equation \eqref{eq:governing}. To this end, we aim to develop a robust and convergent numerical scheme based on a continuous Galerkin-type finite element approach. The efficacy of our model will be demonstrated by solving various configurations with different parameter values and analyzing the resulting numerical predictions.

To properly establish the variational formulation, we introduce the following specialized subspaces of the Sobolev space $\left( W^{1,\,2}(\Omega)\right)^3$. These spaces are crucial for incorporating the specific boundary conditions of our problem:
\begin{subequations}
\begin{align}
V_{\bfzero} &:= \left\{ \bfu \in \left( W^{1,\,2}(\Omega)\right)^3 \colon \; \bfu=\bfzero \quad \text{on} \;\;\Gamma_D\right\}, \label{test_V0}\\
V_{\bfa{g}} &:= \left\{ \bfu \in \left( W^{1,\,2}(\Omega)\right)^3 \colon \; \bfu=\bfa{g} \quad \text{on} \;\; \Gamma_D\right\}. \label{test_Vu0}
\end{align}
\end{subequations}

To derive the weak formulation of the aforementioned BVP, we take the governing equation \eqref{eq:governing}, multiply it by a suitable test function $\bfv \in \widehat{V}_{\bfg}$, and then integrate the resulting expression by parts. Subsequently, applying Green's formula and incorporating the prescribed boundary conditions yields the following continuous weak formulation:

\begin{cwf}
Find $\bfu \in {V}_{\bfg}$, such that
\begin{equation}\label{eq:weak_formulation}
    a(\bfu, \, \bfv) = l(\bfv), \; \forall\, \bfv \in {V}_{\bfzero},
\end{equation}
where the bilinear form $a(\bfu, \, \bfv)$ and the linear form $l(\bfv)$ are, respectively, defined by
\begin{subequations}\label{def:A-L}
\begin{align}
    a(\bfu, \, \bfv) &= \int_{\Omega} \left[  \frac{\mathbb{E}[\bfeps]}{1 + \beta \, \tr(\bfeps)}   \right] \colon \boldsymbol{\varepsilon}( {\bfv}) \; d\bfx\, , \label{eq42a}\\
    l (\bfv) &= \int_{\Omega} \bff \cdot \, \bfv \; d\bfx + \int_{\Gamma_N} \bfg \cdot \bfv \; ds.
\end{align}
\end{subequations}
\end{cwf}

It is crucial to recognize that the continuous weak formulation presented above exhibits nonlinearity. Consequently, to facilitate numerical simulations, particularly for investigating stress-strain concentration in the vicinity of crack-tips within a 3-D body, it becomes necessary to transform this into a sequence of linear problems. We achieve this by constructing a \textit{Picard's type iterative algorithm} at the continuous level. The subsequent linear problems arising from this iteration will then be discretized using a stable Galerkin-type finite element method. This leads to the following iterative algorithm at the continuous level:

\subsubsection*{Picard's Iterative Algorithm} 
For $n=0, \, 1, \, 2, \, \ldots$, and given an initial guess $\bfu^{0} \in {V}_{\bfg}$, find $\bfu^{n+1} \in {V}_{\bfg}$, such that
\begin{equation}\label{eq_pia}
    a(\bfu^{n}; \; \bfu^{n+1}, \, \bfv) = l(\bfv), \; \forall\, \bfv \in {V}_{\bfzero},
\end{equation}
where the terms are defined as:
\begin{subequations}\label{eq_pia_1}
\begin{align}
    a(\bfu^{n}; \; \bfu^{n+1}, \, \bfv) &= \int_{\Omega} \left[ \frac{\mathbb{E}[\bfeps(\bfu^{n+1})]}{1 + \beta \, \tr(\bfeps(\bfu^n))} \right] \colon \boldsymbol{\varepsilon}( {\bfv}) \; d\bfx\, , \label{eq42b}\\ 
    l (\bfv) &= \int_{\Omega} \bff \cdot \, \bfv \; d\bfx + \int_{\Gamma_N} \bfg \cdot \bfv \; ds.
\end{align}
\end{subequations}

\subsection{Finite element space and discrete problem}
\label{sec:fe_discrete_problem}

This section outlines the construction of the {discrete finite element problem}, a numerical counterpart to the continuous weak formulation given by equations \eqref{eq_pia} and \eqref{eq_pia_1}. Our material domain, $\Omega$, is characterized as a 3-D non-convex polyhedral region. We presuppose the existence of a mesh, $\mathcal{T}_h$, where $h>0$ denotes the mesh size. This mesh is assumed to be either quasi-uniform or specifically refined \textit{a priori} to ensure accurate approximation.

The discretization of the domain $\Omega$ strictly adheres to the principles of {conforming and shape-regular finite elements}, as established by Ciarlet \cite{ciarlet2002finite}. This implies that for any two distinct elements $K_1, K_2 \in \mathcal{T}_h$, their intersection $\overline{K}_1 \cap \overline{K}_2$ can only be a null set, a shared vertex, a common edge, or one element wholly contained within another (e.g., $\overline{K}_1 \subset \overline{K}_2$). Furthermore, the union of all elements in the mesh precisely covers the entire domain, i.e., $\bigcup\limits_{K \in \mathcal{T}_h} \overline{K} = \overline{\Omega}$.

For approximating the primary unknown, the displacement field $\bfu$, we introduce a specific functional space,
\begin{equation}
S_h = \left\{ \bfu_h \in \left( C(\overline{\Omega})\right)^3 \colon \left. \bfu_h\right|_K \in \mathbb{Q}_k^d, \; \forall K \in \mathcal{T}_h \right\},
\end{equation}
where $\mathbb{Q}_k^d$ represents a tensor-product polynomial space of order up to $k$ defined over a reference cell $\widehat{K}$. The final discrete approximation space, $\widehat{V}_h$, is then defined as
\begin{equation}\label{app-spaces}
\widehat{V}_h = S_h \, \cap \, V_{\bfg}.
\end{equation}
With these definitions in place, the {discrete finite element problem} is formulated as an iterative scheme:

\paragraph{Discrete Finite Element Problem} 
Given the material parameter $\beta$, the Dirichlet boundary data $\bfu^{0}_h \in \widehat{V}_h$, and the solution from the $n^{th}$ iteration, $\bfu^n_h \in \widehat{V}_h$ (for $n=0, 1, 2, \ldots$), the goal is to find the next iterative solution $\bfu^{n+1}_h \in \widehat{V}_h$ such that:
\begin{equation}\label{discrete-wf}
    a(\bfu_h^n; \, \bfu^{n+1}_h,\, \bfv_h) = l(\bfv_h), \quad \forall\, \bfv_h \in \widehat{V}_h,
\end{equation}
where the bilinear form $a(\cdot\,;\,\cdot\,, \cdot)$ and the linear form $l(\cdot)$ are precisely defined as:
\begin{subequations}\label{A-L-Def}
\begin{align}
a(\bfu_h^n; \, \bfu^{n+1}_h,\, \bfv_h) &=\int_{\Omega} \left[ \frac{\mathbb{E}[\bfeps(\bfu_h^{n+1})]}{1 + \beta \, \tr(\bfeps(\bfu_h^n))} \right] \colon \boldsymbol{\varepsilon}( {\bfv_h}) \; d\bfx\, , \label{disc_A} \\
l(\bfv_h) &= \int_{\Omega} \bff \cdot \, \bfv_h \; d\bfx + \int_{\Gamma_N} \bfg \cdot \bfv_h \; ds. \label{disc_L}
\end{align}
\end{subequations}

This iterative approach constitutes the comprehensive discrete finite element computational procedure. Our overarching aim in this work preserves the accurate characterization of crack-tip fields within a 3-D porous elastic solid.

\begin{enumerate}
    \item {Initialization:}
    \begin{itemize}
        \item Set the iteration counter $n=0$.
        \item Choose an initial guess for the displacement field, $\bfu^0_h \in \widehat{V}_h$. This could be an assumed zero displacement ($\bfu^0_h = \boldsymbol{0}$) or any other suitable initial approximation.
        \item Set a convergence tolerance $\text{TOL} > 0$.
        \item Set a maximum number of iterations $N_{\text{max}}$.
    \end{itemize}

    \item {Iteration Loop:} While $n < N_{\text{max}}$ and convergence is not achieved, perform the following steps:
    \begin{enumerate}
        \item {Formulate the Linearized Problem:} Using the current approximation $\bfu^n_h$, set up the linear system for the next iteration $\bfu^{n+1}_h$. This involves computing the bilinear form $a(\bfu_h^n; \, \cdot \, , \, \cdot)$ and the linear form $l(\cdot)$:
        \begin{align*}
            a(\bfu_h^n; \, \bfu^{n+1}_h,\, \bfv_h) &= \int_{\Omega} \left[ \frac{\mathbb{E}[\bfeps(\bfu_h^{n+1})]}{1 + \beta \, \tr(\bfeps(\bfu_h^n))} \right] \colon \bfeps( {\bfv_h}) \; d\bfx \\
            l(\bfv_h) &= \int_{\Omega} \bff \cdot \, \bfv_h \; d\bfx + \int_{\Gamma_N} \bfg \cdot \bfv_h \; ds.
        \end{align*}
        Note that the term $1 + \beta \, \tr(\bfeps(\bfu_h^n))$ is known from the previous iteration, effectively linearizing the problem with respect to $\bfu^{n+1}_h$.

        \item {Solve the Linear System:} Find $\bfu^{n+1}_h \in \widehat{V}_h$ by solving the discrete linear system:
        \begin{equation}
            a(\bfu_h^n; \, \bfu^{n+1}_h,\, \bfv_h) = l(\bfv_h), \quad \forall\, \bfv_h \in \widehat{V}_h.
        \end{equation}
        This typically translates into solving a matrix system of the form $\mathbf{K}(\bfu_h^n) \mathbf{U}^{n+1} = \mathbf{F}$, where $\mathbf{K}$ is the stiffness matrix (dependent on $\bfu_h^n$), and $\mathbf{U}^{n+1}$ is the vector of nodal degrees of freedom for $\bfu^{n+1}_h$.

        \item {Check for Convergence:} Evaluate a suitable convergence criterion. For instance, check the relative difference between successive solutions:
        \begin{equation*}
            \frac{\|\bfu^{n+1}_h - \bfu^n_h\|}{\|\bfu^{n+1}_h\|} < \text{TOL}
        \end{equation*}
        or the residual of the nonlinear system.

        \item {Update Iteration:} Increment the iteration counter: $n \leftarrow n+1$.
    \end{enumerate}

    \item {Output:} The converged solution $\bfu_h^{N}$, where $N$ is the final iteration number.
\end{enumerate}

\subsubsection{Postprocessing and key mechanical variables}
\label{sec:postprocessing} 

In this study, a critical aspect of characterizing the intricate behavior of the novel porous elastic solid involves performing several important postprocessing steps. These computations are essential for understanding the {crack-tip fields} and extracting meaningful insights from our simulations. Specifically, the following key mechanical variables are meticulously calculated using the final converged displacement solution, $\bfu_h$:
\begin{subequations}
\begin{align}
\bfT_h &= \frac{\overline{c}_1}{2} \, \left( \nabla \bfu_h + (\nabla \bfu_h)^T \right) + \overline{c}_2 \, (\nabla \cdot \bfu_h) \, \bfI \label{eq:stress_tensor_post}\\
\bfeps_h &= \frac{(1+\nu)(1+ \beta \, \nabla \cdot \bfu_h)}{E} \, \bfT_h - \frac{\nu(1+ \beta \, \nabla \cdot \bfu_h )}{E} \, \tr(\bfT_h) \, \bfI. \label{eq:strain_tensor_post}
\end{align}
\end{subequations}
Beyond these fundamental quantities, another crucial variable for understanding {crack-tip behavior} is the \textit{strain energy density}. This quantity, which symbolizes the material's strength and directly quantifies crack-tip evolution, is derived by taking the tensor inner product of the computed stress ($\bfT_h$) and strain ($\bfeps_h$) tensors.

This work represents a fundamental effort to comprehensively characterize the 3-D behavior of a porous elastic solid whose material properties are explicitly dependent on its density. This is a significant contribution to the field, as it provides a robust framework for analyzing such complex material responses.

\section{Results and discussion}\label{rd}
\begin{figure}[H]
	\centering
	\subfloat[\footnotesize{The star-shaped crack and the Cartesian coordinate system with its three axes in the directions of the three perpendicular edges (length, width, and height).}]{\label{star}\includegraphics[width=0.4\textwidth]{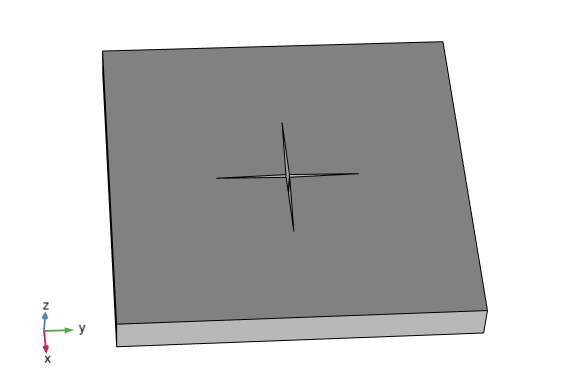}} $\qquad$
	\subfloat[\footnotesize{The tetrahedral mesh of the geometry. Mesh is refined near the tips of the star for more accurate computational results. }]{\label{starmesh}\includegraphics[width=0.45\textwidth]{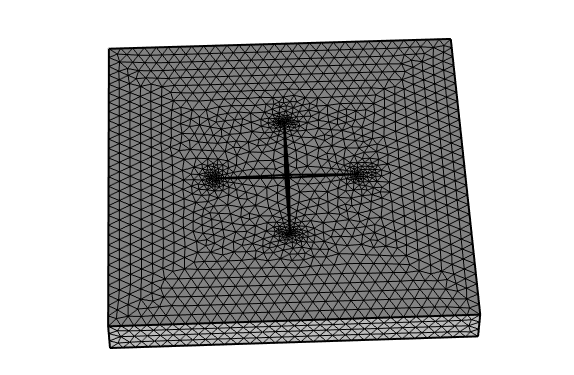}} 
	\caption{{\footnotesize{The star-shaped crack embedded in a square plate and tetrahedral meshes of the geometry.}}}
	\label{}
\end{figure}

The star-shaped crack is embedded in a square plate (Fig. \ref{star}). The rectangular coordinate system is oriented with its three axes in the direction of the three mutually perpendicular edges for length, width, and height. The side of the square is 10 cm long. The thickness of the plate is 1 cm. Each star-tip angle is about $3^{\text{o}}$.  The tetrahedral mesh of the geometry is shown in Fig. \ref{starmesh}. A total number of 49168 tetrahedrons is constructed in the mesh with a total of 7408 triangular meshes for the surface, and  a total number of 480 line-segment mesh on the edge.

%

\begin{figure}[H]
		\centering
		\subfloat[\footnotesize{The displacement boundary conditions on the left and right surfaces with the displacement parallel to the $y$-direction.}]{\label{boundary}\includegraphics[width=0.4\textwidth]{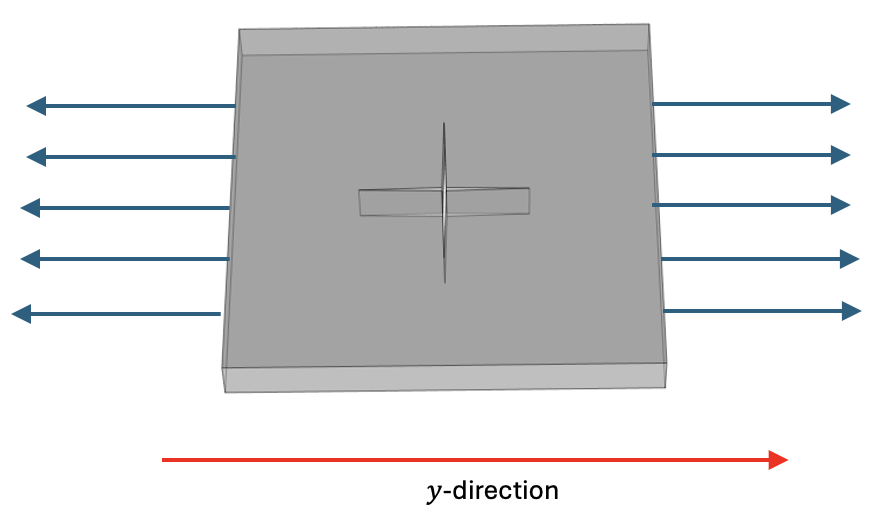}} $\qquad$
		\subfloat[\footnotesize{The $r$-direction along which mechanical analysis is performed. The origin of the $r$-direction is in the middle of the tip and the $r$-direction is  parallel to the $x$-direction. }]{\label{direction}\includegraphics[width=0.45\textwidth]{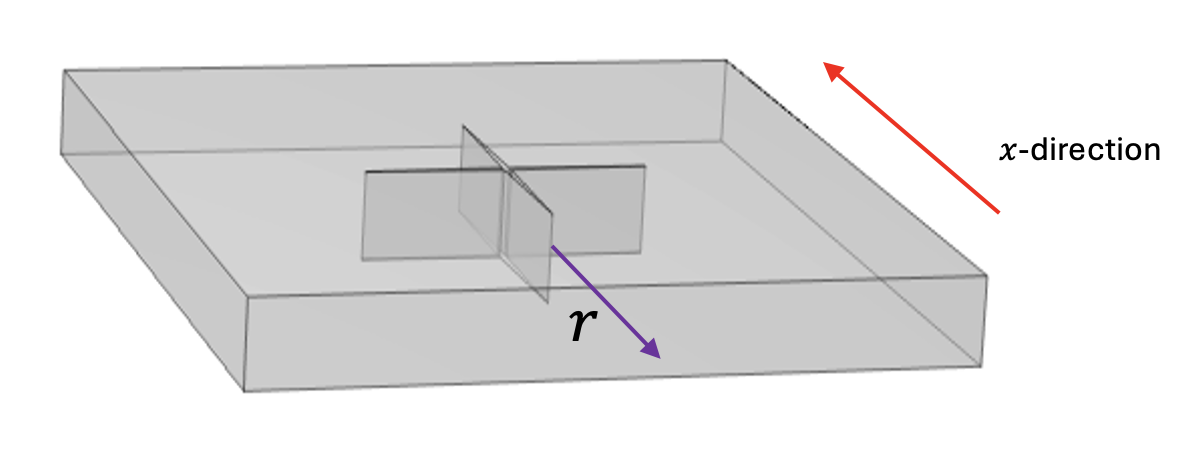}} 
			\caption{{\footnotesize{The displacement boundary condition and the $r$-direction along which mechanical analysis is performed.}}}
		\label{georadi}
	\end{figure}

The displacement boundary condition is imposed on the two side surfaces as shown in Fig. \ref{boundary}, and the displacements on these two surfaces are parallel to the $y$-direction. On the left surface, the displacement is $-1$ mm, and on the right surface, the displacement is 1 mm. The displacements for the boundary conditions are symmetric for force equilibrium to the plate. Mechanical analysis for $T_{22}$ and $\epsilon_{22}$ is performed along the $r$-axis as shown in Fig. \ref{direction}. The origin of the $r$-axis is in the middle of the crack tip, and the $r$-axis is parallel to the $x$-direction.

\begin{figure}[H]
	\centering
	\subfloat[\footnotesize{$T_{22}$  (unit: $10^4$ Pa).}]{\label{t22}\includegraphics[width=0.37\textwidth]{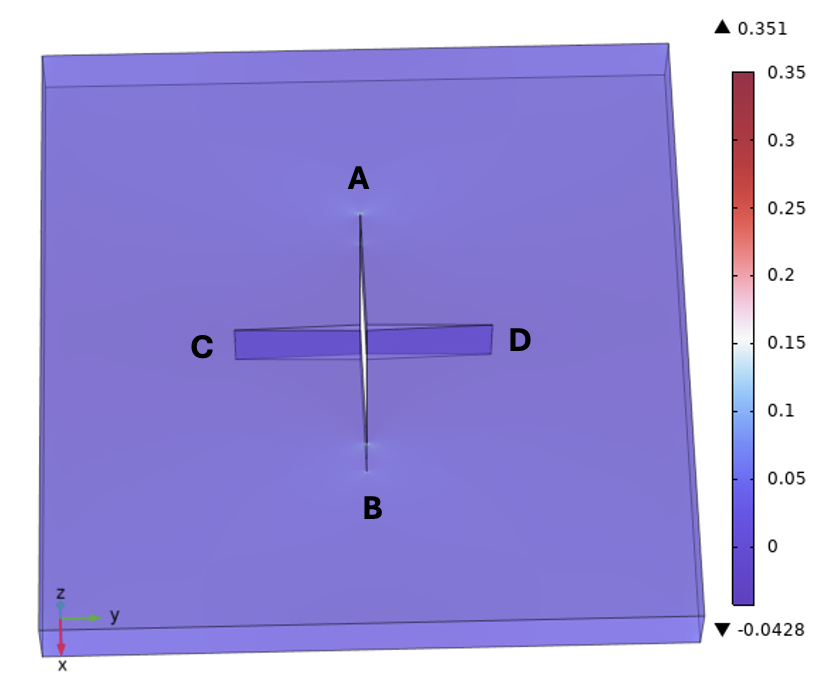}} $\qquad$
	\subfloat[\footnotesize{ $\epsilon_{22}$.}]{\label{s22}\includegraphics[width=0.35\textwidth]{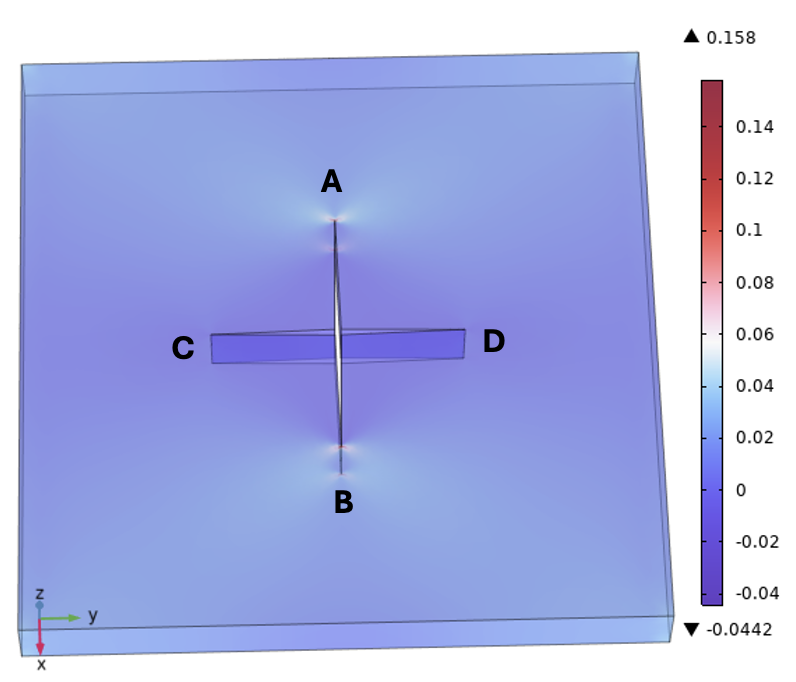}} 
	\caption{{\footnotesize{$T_{22}$ and $\epsilon_{22}$ distributions for $\beta=-2$ in the  plate  in a transparent-view mode. Stress and strain concentration can be viewed at the tips A and B due to the boundary conditions along the $y$-axis to open the cracks. Tips C and D show no stress and strain concentration due to the effect of the boundary conditions  to close the cracks rather than to open them.}}}
	\label{concente}
\end{figure}

Computational results for all $\beta$ values show concentration of stress and strain near the two crack tips A and B with their pertinent  cracks in the $x$ direction (Fig. \ref{concente}). We illustrate the distribution for stress $T_{22}$ and strain $\epsilon_{22}$  for $\beta=-2$ on the plate  as an example. Concentration represented by brighter spots  can be viewed near the crack tips A and B. Fig. \ref{t22} shows a stress concentration weaker than Fig. \ref{s22} for strain concentration. Crack tips C and D show no concentration because the boundary conditions tend to close the cracks rather than to open them. Through the transparent-view mode of the plate, over the crack tips A and B, concentrations are formed along the whole tip from the top to the bottom, but locations in the top and bottom surfaces show the strongest concentration.

\begin{figure}[H]
	\centering
	\subfloat[\footnotesize{$T_{22}$ (unit: $10^4$ Pa).}]{\label{stress22neg}\includegraphics[width=0.45\textwidth]{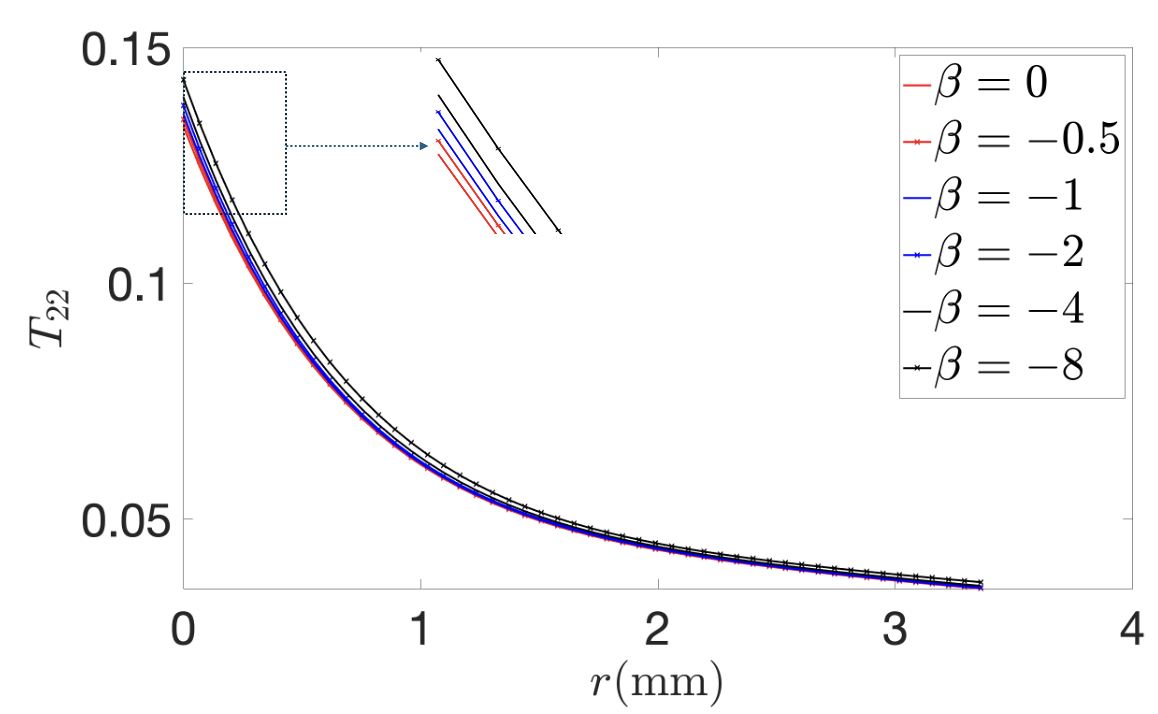}} $\qquad$
	\subfloat[\footnotesize{ $\epsilon_{22}$.}]{\label{strain22gena}\includegraphics[width=0.45\textwidth]{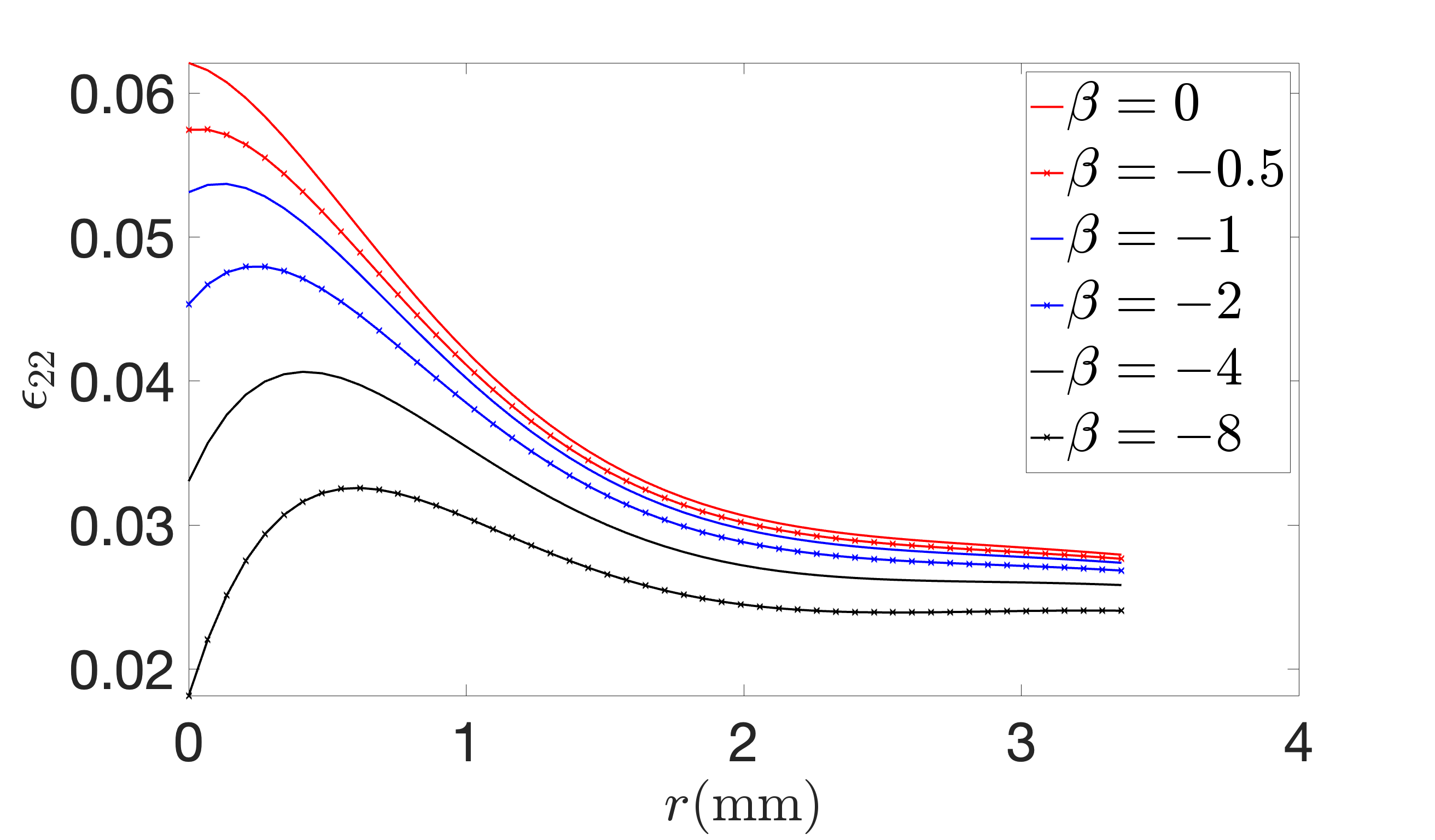}} 
	\caption{{\footnotesize{$T_{22}$ and $\epsilon_{22}$ distributions for negative $\beta$ values $(0,\, -0.5,\, -1,\, -2,\, -4,\, -8)$ along the $r$ line in Fig. \ref{direction}. The curves for $T_{22}$ in the left panel is dense near the crack tip with $r=0$. The curves near the tip is magnified and illustrated separately for better view. Near the crack tip at $r=0$, $\beta$ further away from 0 increases $T_{22}$ but decreases $\epsilon_{22}$.}}}
	\label{distriis}
\end{figure}

Distribution of $T_{22}$ and $\epsilon_{22}$ for mainly negative $\beta$ values $(0,\, -0.5,\, -1,\, -2,\, -4,\, -8)$ along the $r$-line (Fig. \ref{direction}) is demonstrated in Fig. \ref{distriis}.  The result for $\beta=0$ is provided as an outcome from the linear elasticity fracture model for comparison with results for non-zero $\beta$ values. The left panel Fig. \ref{stress22neg} shows distribution for $T_{22}$. As $\beta$ decreases from 0 to $-8$, $T_{22}$ value is increasing near the crack tip ($r=0$). For any curve with a fixed $\beta$ value, $T_{22}$ is decreasing as $r$ increases. The curves for $T_{22}$ with different $\beta$ values are layered and do not intersect each other.  As $\beta$ is decreasing from 0 to $-8$, $\epsilon_{22}$  is decreasing near the crack tip ($r=0$) shown in Fig. \ref{strain22gena}. Different from $T_{22}$, for a fixed $\beta$, $\epsilon_{22}$ is not monotonic. $\epsilon_{22}$ first increases and then decreases, forming a bump in the curve. Such bump is more pronounced when $\beta$ is further away from 0. The curves for $\epsilon_{22}$ with different $\beta$ values are also layered with no intersection.   Different $\beta$ effect on $T_{22}$ and $\epsilon_{22}$ demonstrates the effect of the nonlinear elasticity model in controlling the amount of stress and strain differently.

Table \ref{percenega} demonstrates the values of $T_{22}$ and $\epsilon_{22}$  near the crack tip $r=0$ for different negative $\beta$ values. In particular, the table shows percentage of increase for $T_{22}$ compared to  the corresponding one at $\beta=0$ (the third row). All the relative increase for $T_{22}$ is positive, representing positive increase of $T_{22}$ as $\beta$ decreases. Comparative increase for $\epsilon_{22}$ is shown in the fifth row with all negative values, representing decrease of $\epsilon_{22}$ as $\beta$ decreases. Only considering the absolute value of the relative change, $\epsilon_{22}$ is changing more greatly than $T_{22}$ when compared to their corresponding values at $r=0$.
\begin{table}[H]
	\begin{center}
		\begin{tabular}{  lcccccc}
			\hline
		$\beta$:	&0  & $-0.5$  &$-1$ &$-2$&$-4$&$-8$   \\
		$T_{22}|_\beta$ (Unit: Pa): &0.130&0.132&0.133&0.135&0.138&0.143\\
		$\frac{T_{22}|_{\beta}-T_{22}|_{\beta=0}}{T_{22}|_{\beta=0}}$: &0\%&1.54\%&2.31\%&3.85\%&6.15\%&10.00\%\\
			$\epsilon_{22}|_\beta$: &0.062&0.058&0.054&0.047&0.035&0.021\\
	$	\frac{\epsilon_{22}|_{\beta}-\epsilon_{22}|_{\beta=0}}{\epsilon_{22}|_{\beta=0}} $: &0\% &  $-6.45\%$&  $-12.90\%$&  $-24.19\%$&  $-43.55\%$&  $-66.13\%$\\
				\hline
		\end{tabular}
		\caption{\footnotesize{$T_{22}$ and $\epsilon_{22}$ values near $r=0$ for different negative $\beta$ values. The table also shows percentage of increase for $T_{22}$ and $\epsilon_{22}$ at $\beta\neq0$ compared to the corresponding ones at $\beta=0$. }}
		\label{percenega}
	\end{center}
\end{table}

\begin{figure}[H]
	\centering
	\subfloat[\footnotesize{$T_{22}$ (unit: $10^4$ Pa).}]{\label{stress22_posi}\includegraphics[width=0.45\textwidth]{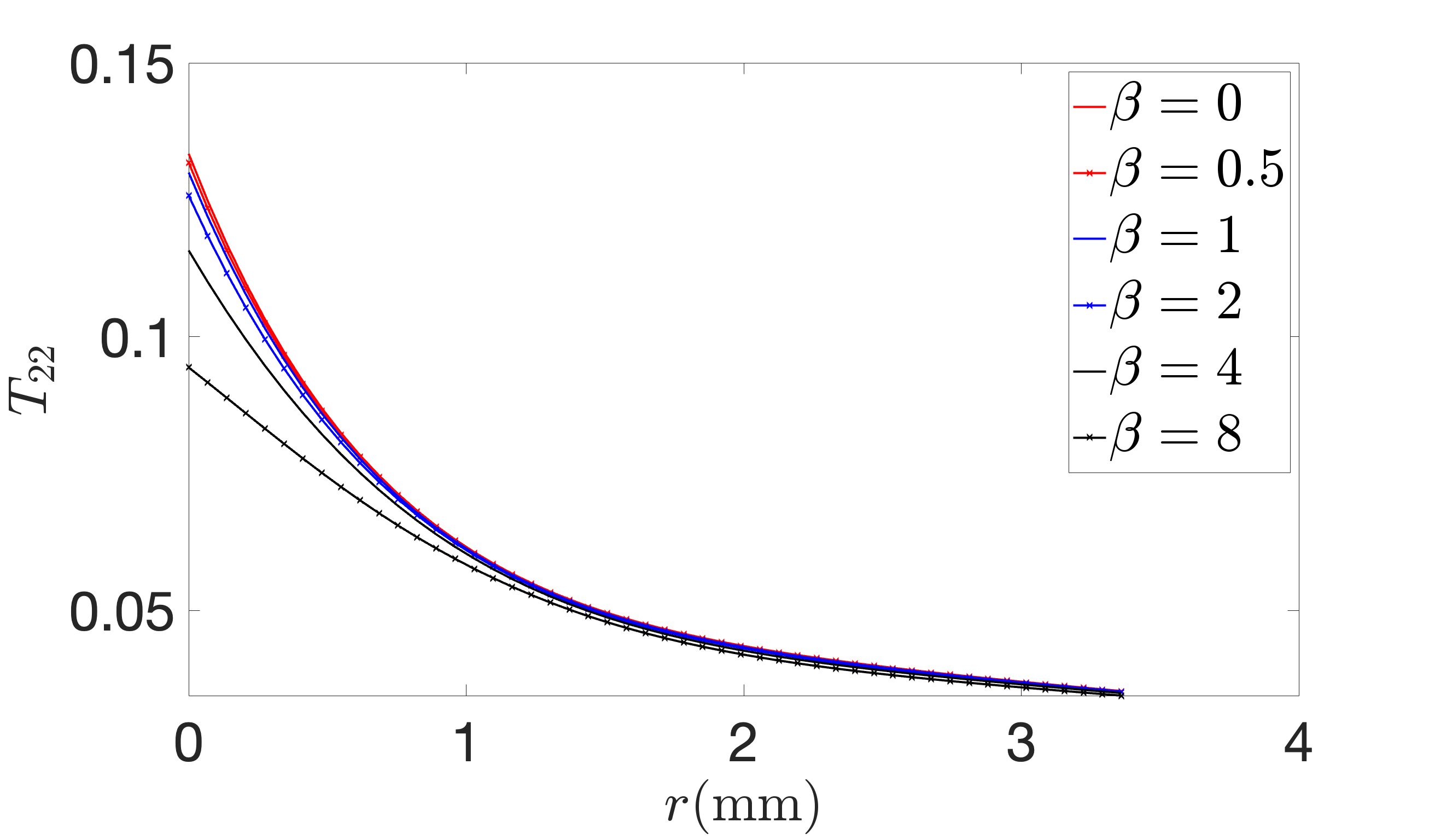}} $\qquad$
	\subfloat[\footnotesize{ $\epsilon_{22}$.}]{\label{strain22_positive}\includegraphics[width=0.45\textwidth]{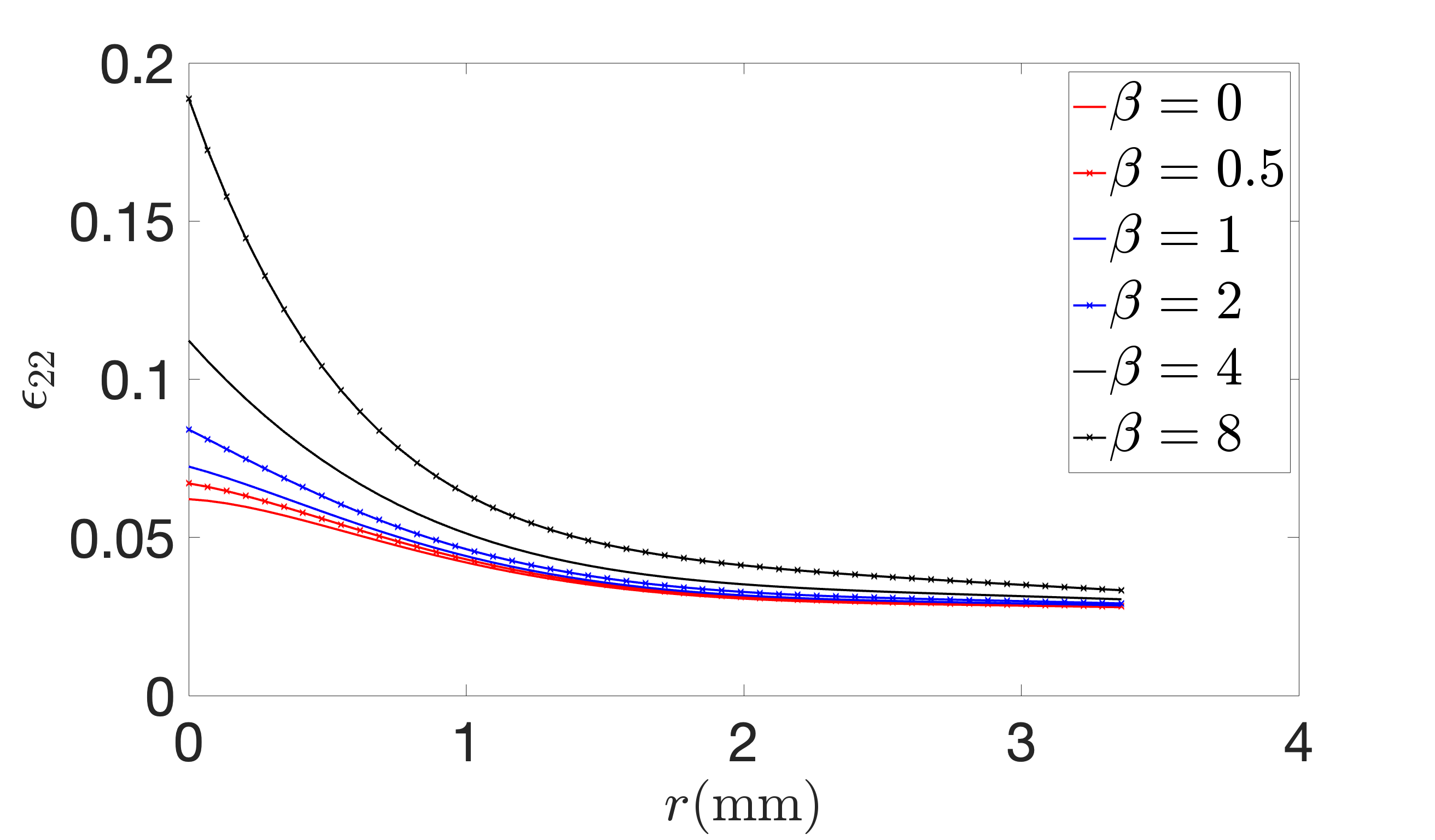}} 
	\caption{{\footnotesize{$T_{22}$ and $\epsilon_{22}$ distributions for positive $\beta$ values $(0,\, 0.5,\, 1,\, 2,\, 4,\, 8)$ along the $r$ line in Fig. \ref{direction}.  Near the crack tip at $r=0$, $\beta$ further away from 0 decreases $T_{22}$ but increases $\epsilon_{22}$, a trend opposite to results in Fig. \ref{distriis} for negative $\beta$ values.}}}
	\label{distriisPosiBet}
\end{figure}

Distribution of $T_{22}$ and $\epsilon_{22}$ for positive $\beta$ values $(0,\, 0.5,\, 1,\, 2,\, 4,\, 8)$ along the $r$-line (Fig. \ref{direction}) is demonstrated in Fig. \ref{distriisPosiBet}. The left panel Fig. \ref{stress22_posi} shows the distribution of $T_{22}$. Near the crack tip at $r=0$, when $\beta$ increases, $T_{22}$ decreases. On the contrary, when $\beta$ increases, $\epsilon_{22}$ increases (Fig. \ref{strain22_positive}). 
Similar to exhibition in Table \ref{percenega} for negative $\beta$ values, Table \ref{percepositive} shows the results for positive $\beta$ values. The third row of the table shows negative values for relative increase of $T_{22}$, representing decease of $T_{22}$ as $\beta$ increases. The fifth row shows positive values for relative increase of $\epsilon_{22}$, representing increase of $\epsilon_{22}$ as $\beta$ increases. The absolute value  $|\frac{T_{22}|_{\beta}-T_{22}|_{\beta=0}}{T_{22}|_{\beta=0}}|$ for each positive $\beta$ is larger than that for the corresponding  negative $\beta$. The same is true for $	|\frac{\epsilon_{22}|_{\beta}-\epsilon_{22}|_{\beta=0}}{\epsilon_{22}|_{\beta=0}}| $. Such outcomes show that a positive $\beta$ value can cause more drastic change than the corresponding negative $\beta$ in the nonlinear elasticity fracture model.
\begin{table}[H]
	\begin{center}
		\begin{tabular}{  lcccccc}
			\hline
			$\beta$:	&0  & 0.5  &1 &2&4&8   \\
			$T_{22}|_\beta$ (Unit: Pa): &0.130&0.128&0.126&0.121&0.111&0.091\\
			$\frac{T_{22}|_{\beta}-T_{22}|_{\beta=0}}{T_{22}|_{\beta=0}}$: &0\%& $-1.54\%$&   $-3.08\%$&   $-6.92\%$&  $-14.62\%$&  $-30.00\%$\\
			$\epsilon_{22}|_\beta$: &0.062&0.067&0.071&0.082&0.108&0.178\\
			$	\frac{\epsilon_{22}|_{\beta}-\epsilon_{22}|_{\beta=0}}{\epsilon_{22}|_{\beta=0}} $: &0\%&    8.06\%&   14.52\%&   32.26\%&   74.19\%&  187.10\% \\
			\hline
		\end{tabular}
		\caption{\footnotesize{$T_{22}$ and $\epsilon_{22}$ values near $r=0$ for different positive $\beta$ values. The table also shows percentage of increase for $T_{22}$ and $\epsilon_{22}$ at $\beta\neq0$ compared to the corresponding ones at $\beta=0$. }}
		\label{percepositive}
	\end{center}
\end{table}

\begin{figure}[H]
	\centering
	\subfloat[\footnotesize{Energy for negative $\beta$ values.}]{\label{energynega}\includegraphics[width=0.45\textwidth]{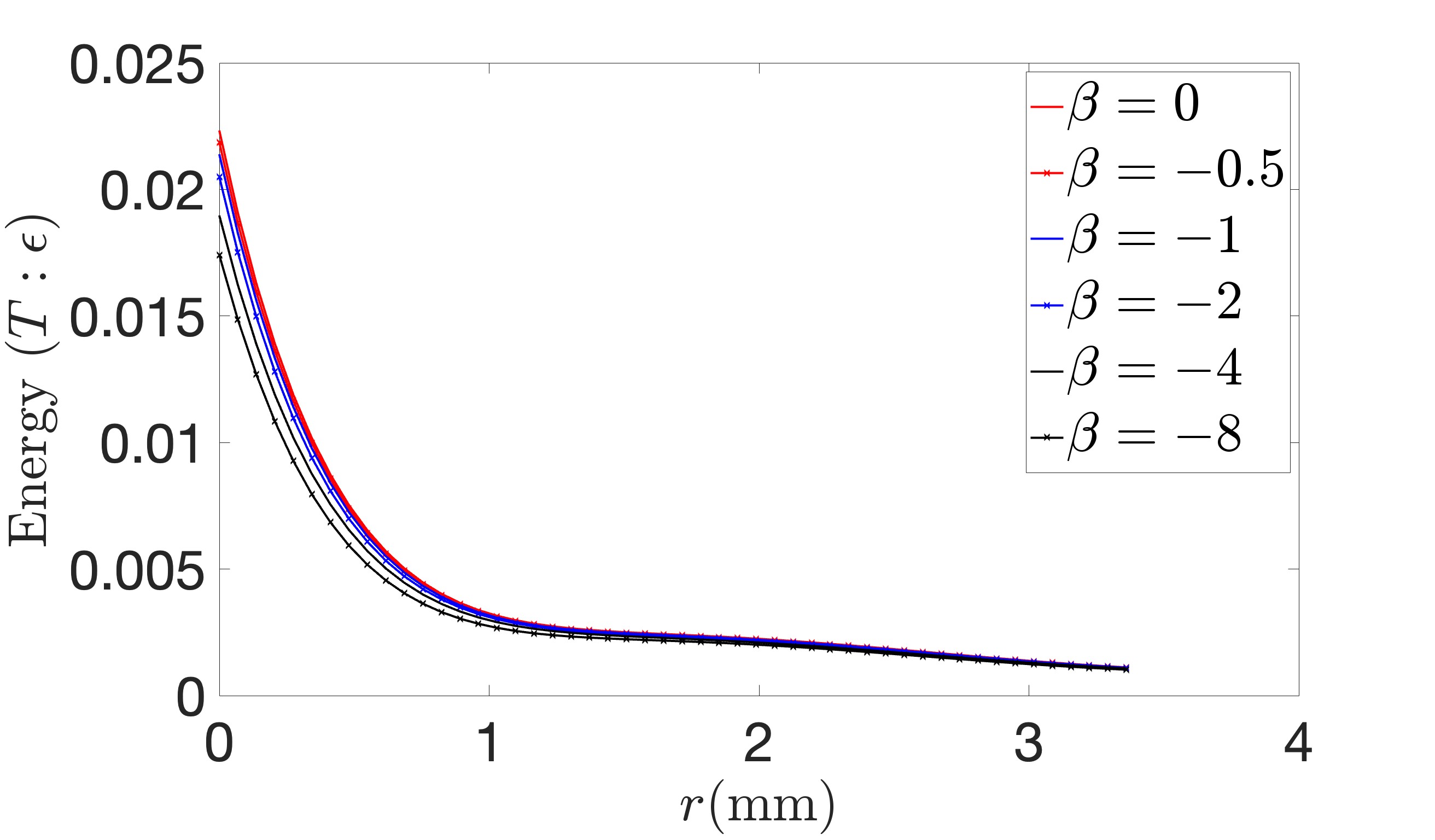}} $\qquad$
	\subfloat[\footnotesize{ Energy for positive $\beta$ values.}]{\label{energyposi}\includegraphics[width=0.45\textwidth]{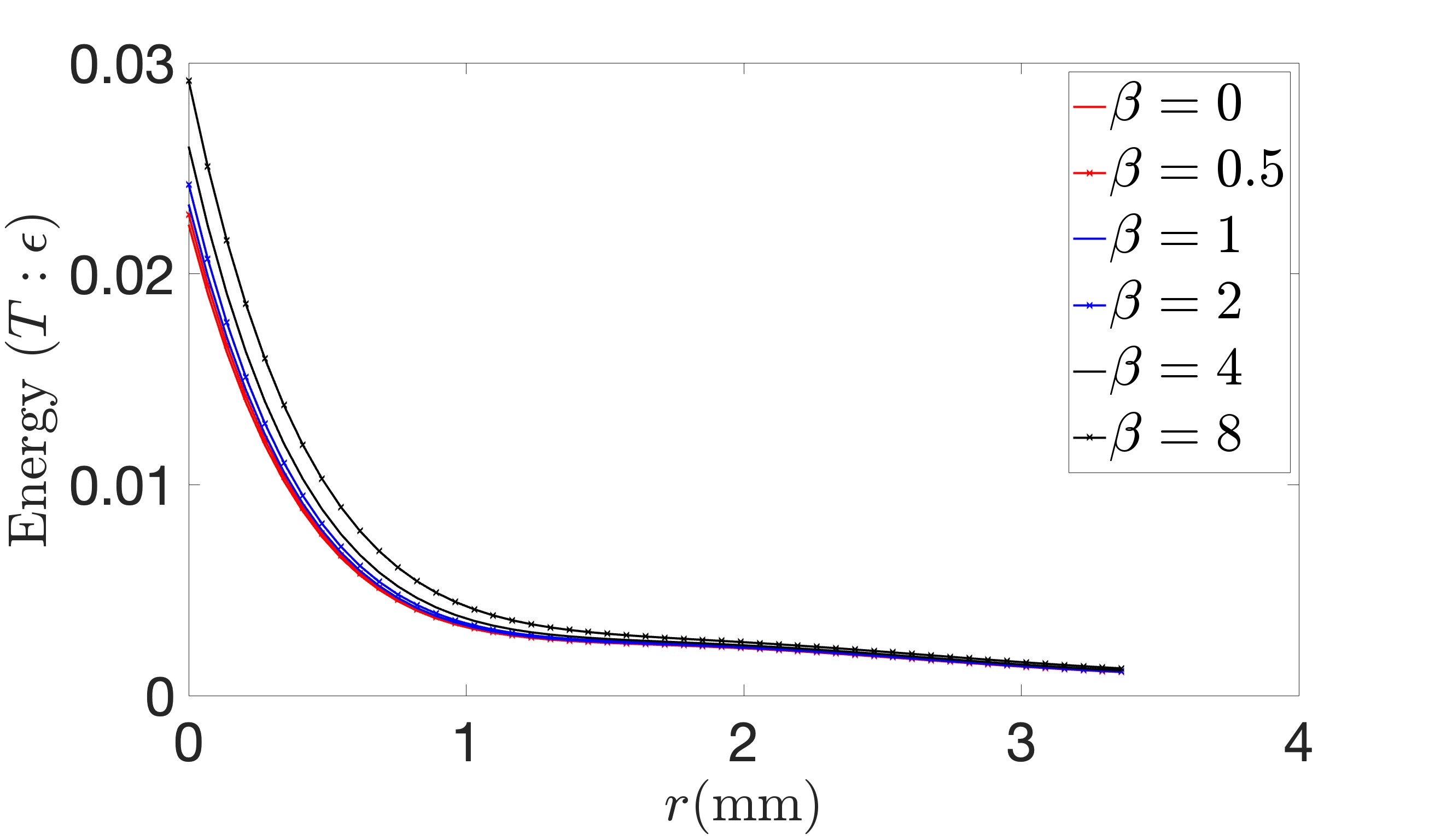}} 
	\caption{{\footnotesize{Energy distributions along the $r$-line for negative and positive $\beta$ values (energy unit: $10^4$mm$^{1/2}$Pa). For negative $\beta$ values, as $\beta$ is further away from 0, the energy decreases near the crack tip. For positive $\beta$ values, as $\beta$ is further away from 0, the energy increases near the crack tip.}}}
	\label{energies}
\end{figure}

Distribution of energy over the whole plate is similar to $T_{22}$ and $\epsilon_{22}$ in Fig. \ref{concente} with concentration near the crack tips at A and B. For brevity, we skip displaying the energy distribution over the whole plate. Figure \ref{energies} shows distribution of energy along the $r$-line. Compared with energy for $\beta=0$ (energy for the linear elasticity fracture model), negative $\beta$ tends to reduce the energy while positive $\beta$ tends to increase the energy, a trend similar to the results for $\epsilon_{22}$ but opposite to the trend for $T_{22}$.

\begin{figure}[H]
	\centering
	\includegraphics[width=0.6\textwidth]{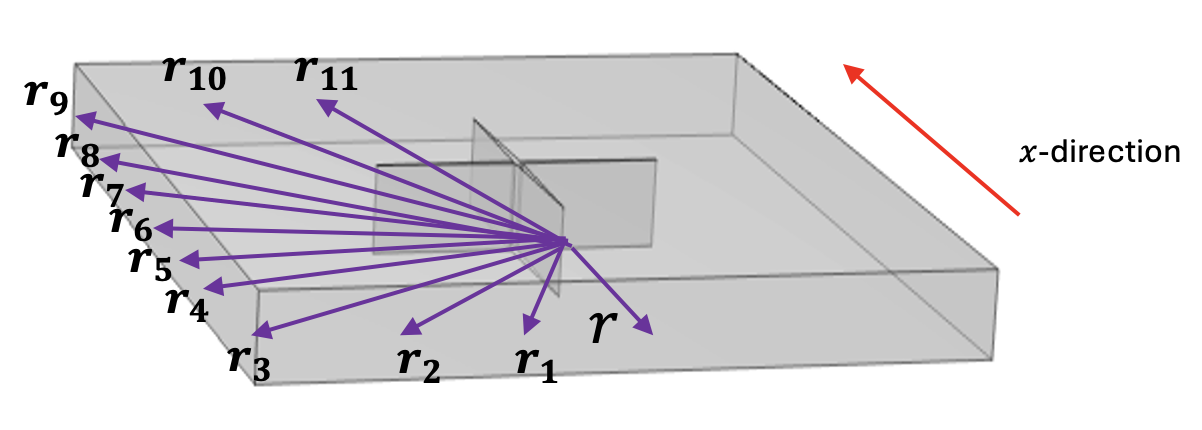}
	\caption{Another 11 directions to study energy distributions along different directions. These directions are denoted, respectively, by $r_1$, $r_2$, $r_3$, $\cdots$, $r_{11}$.}
	\label{direcs}
\end{figure}

The energy around the crack tip may not be evenly distributed in different directions.  In Fig. \ref{direcs}, we draw another 11 directions, denoted, respectively, by $r_1$, $r_2$, $r_3$, $\cdots$, $r_{11}$, to more deeply study the energy distribution. All the 11 directions start from the middle crack tip as for $r$, and ends at  the plate's boundary (middle of the thickness). The other half of the plate without these directions is symmetric in distribution of energy. The energy contour distribution is displayed in Table \ref{energydistss}. The energy values are for the location near the origin tip of each $r_i$ line (not exactly at the tip of each $r_i$) due to numerical processing for discrete meshes. For some neighboring $r_i$'s, the energy values may be identical due to difficulty of differentiation of too dense meshes.  The following are some main results observed from the table: (1) Along the $r$-direction, as $\beta$ increases from $-8$ to 8, the energy also increases; (2) For any fixed $\beta$, energy near the tip of $r_2$ (or $r_3$) is the greatest in the corresponding column; (3) For $\beta=0$, energy  increases and then decreases from $r_1$ to $r_8$, with the energy maximum at $r_2$ (or $r_3$); (4) Energy near the tip of $r_1$, $r_5$, or $r_6$ increases first and then decreases as $\beta$ increases from $-8$ to 8, with its maximum at $\beta=0$; (5) Near the tips of $r_2$, $r_3$, or $r_4$, energy decreases as $\beta$ increases.

\begin{table}[H]
	\small
	\begin{center}
		\begin{tabular}{  lccccccccccc}
			\hline
			$\beta$:	&$-8$&$-4$&$-2$&$-1$&$-0.5$&0&0.5&1&2&4&8   \\
			$r$: &0.0179&0.0192&0.0205&0.0213&0.0217&0.0221      &0.0225&0.0229&0.0237&0.0251&0.0276\\
			$r_1$: &0.0179&0.0192&0.0205&0.0213&0.0217& 0.0276    &0.0251&0.0237&0.0229&0.0225&0.0221\\
			$r_2$: &0.1762&0.1142&0.0957&0.0885&0.0853&0.0824   &0.0797&0.0772&0.0728&0.0658&0.0565\\
			$r_3$:&0.1762&0.1142&0.0957&0.0885&0.0853&0.0824     &0.0797&0.0772&0.0728&0.0658&0.0565\\
			$r_4$: &0.0954&0.0867&0.0796&0.0761&0.0744& 0.0727     &0.0711&0.0695&0.0666&0.0615&0.0539\\
			$r_5$: &0.0410&0.0482&0.0495&0.0496&0.0496& 0.0494     &0.0493&0.0490&0.0485&0.0471&0.0442\\
			$r_6$: &0.0410&0.0482&0.0495&0.0496&0.0496& 0.0494     &0.0493&0.0490&0.0485&0.0471&0.0442\\
			$r_7$: &0.0437&0.0417&0.0410&0.0406&0.0405& 0.0403     &0.0402&0.0401&0.0398&0.0396&0.0401\\
			$r_8$: &0.0437&0.0417&0.0410&0.0406&0.0405& 0.0403     &0.0402&0.0401&0.0398&0.0396&0.0401\\
			$r_9$: &0.0437&0.0417&0.0410&0.0406&0.0405& 0.0403     &0.0402&0.0401&0.0398&0.0396&0.0401\\
			$r_{10}$: &0.0464&0.0439&0.0426&0.0418&0.0414& 0.0410    &0.0407&0.0403&0.0397&0.0387&0.0383\\
			$r_{11}$: &0.0464&0.0439&0.0426&0.0418&0.0414& 0.0410   &0.0407&0.0403&0.0397&0.0387&0.0383\\
			\hline
		\end{tabular}
		\caption{\footnotesize{ Energy (unit: $10^4$mm$^{1/2}$Pa) near the origin  of the 12 directions  ( $r=0$ and $r_i=0$, for $i$=1, 2, $\cdots$, 11) for both negative and positive $\beta$ values for a contour analysis of the energy distribution around the crack tip.}}
		\label{energydistss}
	\end{center}
\end{table}

\section{{ Conclusion}}\label{sec:NumExp}
The past few years have witnessed a surge in theoretical investigations aimed at developing algebraically nonlinear models. These models are designed to accurately describe both geometrically linear and nonlinear elastic solids, as evidenced by foundational works such as \cite{rajagopal2003,rajagopal2007elasticity}. The overarching objective of these rigorous studies is to formulate mathematically well-posed constitutive theories that govern the bulk material behavior. Crucially, these advanced models seek to predict {bounded crack-tip strain fields}, thereby circumventing the need for \textit{ad hoc} modeling interventions, such as cohesive or process zones, which are often employed to rectify the unphysical singular predictions arising from classical linear elasticity. Perhaps the most profound implication of the modeling paradigm established in \cite{rajagopal2011modeling,gou2015modeling,gou2023computational,gou2023finite,Mallikarjunaiah2015,MalliPhD2015} lies in its potential to revolutionize fracture theories. This new framework paves the way for developing more robust and physically consistent approaches applicable to a diverse range of brittle materials.

This paper represents the inaugural phase of a comprehensive research initiative dedicated to establishing both theoretical and computational models for cracks and fractures within a novel class of 3-D porous elastic solids. These materials are particularly complex as their mechanical response is intrinsically linked to their density. Specifically, leveraging a particular subclass of these constitutive relationships, we herein derive a quasilinear elliptic BVP tailored for static analysis within a star-shaped crack. For the numerical solution of this 3-D problem, we employ a {continuous Galerkin-type finite element method}. This approach is seamlessly integrated with a \textit{Picard-type linearization scheme}, a crucial coupling that ensures the problem remains linear at each iterative step. This linearization is highly advantageous, as it guarantees the existence of a unique solution for every iteration. Consequently, we can reliably use a direct solver to accurately determine the nodal values of the displacement field, leading to robust computational results.

Our observations confirm that the finite element discretization of the 3-D BVP exhibits excellent convergence. The Picard iterations, which linearize the problem at each step, consistently reach the defined stopping criterion within a reasonable number of iterations. Specifically, we set the {relative tolerance} for convergence at $0.001$, with a generous {maximum limit of $1000$ iterations} to ensure robustness across various scenarios.

A crucial aspect of our numerical approach is the {initial guess} for the Picard iteration. This is intelligently derived by first solving the BVP under {linearized elasticity} assumptions, effectively setting the material parameter $\beta=0$ in our more comprehensive nonlinear model. This strategy provides a physically informed starting point, contributing to the efficiency of the iterative process. Furthermore, our results demonstrate a significant and expected physical behavior: as the nonlinearity parameter $\beta$ approaches zero, the numerical solution obtained from our proposed nonlinear BVP converges to the corresponding solution predicted by the {classical linear elasticity model}. This provides strong validation for our model and its implementation, affirming its consistency with established theory in the appropriate limits.
A few remarks about the numerical solution and its interpretation are as follows.
\begin{itemize}
\item[1.]  Under tensile loading conditions, our nonlinear model exhibits distinct and insightful behavior depending on the sign of the material parameter $\beta$. For {negative values of $\beta$}, the nonlinear model predicts a significantly {attenuated crack-tip strain growth} compared to the predictions of the classical linear elastic model. This observation highlights a clear {strain-limiting effect} inherent to our proposed model when $\beta$ is negative, suggesting a toughening mechanism at the crack tip. Conversely, a precisely {opposite trend} in crack-tip strain values is observed when $\beta$ takes on {positive values}. In this scenario, the nonlinear effects lead to a different response at the crack vicinity. This contrasting behavior for positive and negative $\beta$ is fundamentally consistent with the model's formulation: the effective stiffness of the material, represented by the term $(1 + \beta \, \tr \, \bfeps)$, is no longer a constant, unlike in the classical linearized elasticity model. Its dynamic variation directly dictates how strain and stress fields evolve.

Regarding {stress concentration}, our analysis reveals that as the negative value of $\beta$ decreases (i.e., becomes more negative), the stress concentration at the crack tip intensifies. This effect aligns with observations from classical models, implying that for smaller (more negative) $\beta$ values, the proposed material model accentuates a {stronger brittle response}. Notably, the largest value of the stress component $\bfT_{22}$ consistently appears directly \textit{before} the crack tip, which is another characteristic feature often observed in classical linear elastic solutions. However, a divergence from classical behavior is noted for {positive $\beta$ values}: in the nonlinear model, stress concentration \textit{decreases} as positive $\beta$ values increase. This indicates a different interplay between material nonlinearity and stress redistribution for different parameter regimes.

\item[2. ]Our analysis reveals that the {strain energy density} also reaches its peak value just ahead of the crack tip. This finding aligns perfectly with predictions from the {classical linear elasticity model}. This consistency is crucial; it means that established local fracture criteria, traditionally applied to linear elastic materials, remain valid and can be directly utilized within the framework of the nonlinear theory presented in this paper to investigate crack evolution problems.

\item[3. ] Furthermore, this framework can be extended to study {quasi-static crack propagation} by incorporating phase-field regularization techniques. Such an approach mirrors methods successfully applied to strain-limiting models, as demonstrated in works like \cite{yoon2021quasi,lee2022finite,manohar2025adaptive,manohar2025convergence,fernando2025textsf,fernando2025xi,ghosh2025computational,ghosh2025finite}. This opens up promising avenues for future research into the dynamic behavior of cracks in these novel porous elastic solids.
\end{itemize}

\section{Acknowledgement.}
The work of SMM  is supported by the National Science Foundation under Grant No. 2316905. 

\bibliographystyle{unsrtnat}
\bibliography{references}

\end{document}